\documentclass[a4paper]{ifacconf}

\usepackage[pdftex]{graphicx}
\graphicspath{./Figs/}

\usepackage{amsmath, amssymb, amsfonts}
\usepackage{xcolor}
\usepackage[round]{natbib}

\usepackage{enumitem}

\renewcommand{\qed}{\hfill\ensuremath{\blacksquare}}

\usepackage{soul}

\usepackage[utf8]{inputenc}

\usepackage{pgfplots}
\usepackage{grffile}
\pgfplotsset{compat=newest}
\usetikzlibrary{plotmarks}
\usetikzlibrary{arrows.meta}
\usepgfplotslibrary{patchplots}

\newcommand\oprocendsymbol{\hbox{\small $\blacksquare$}}
\newcommand\oprocend{\relax\ifmmode\else\unskip\hfill\fi\oprocendsymbol}

\newtheorem{theorem}{Theorem}
\newtheorem{lemma}{Lemma}
\newtheorem{definition}{Definition}
\newtheorem{assumption}{Assumption}
\newtheorem{proposition}{Proposition}
\newtheorem{corollary}{Corollary}
\newtheorem{remark}{Remark}

\DeclareMathOperator{\cl}{cl}

\DeclareMathOperator{\diag}{diag}

\newcommand{\calX}{\mathcal{X}}

\newcommand{\calZ}{\mathcal{Z}}
\newcommand{\calM}{\mathcal{M}}
\newcommand{\calL}{\mathcal{L}}
\newcommand{\bbR}{\mathbb{R}}

\newcommand{\bbI}{\mathbb{I}}

\newcommand{\calS}{\mathcal{S}}
\newcommand{\calF}{\mathcal{F}}

\begin{document}
\begin{frontmatter}

  \title{Limit Behavior and the Role of Augmentation in Projected Saddle Flows for Convex Optimization}

  \author{
    \hspace{-1.25cm}
    Adrian Hauswirth$^1$, Lukas Ortmann$^1$, Saverio Bolognani, Florian D\"orfler}

  \address{Department of Information Technology and Electrical Engineering, ETH Zurich,
  8092 Zurich, Switzerland
  Email:{\tt \{hadrian,ortmannl,bsaverio,dorfler\}@ethz.ch}}

  \thanks{A. Hauswirth and L. Ortmann contributed equally to this paper.}
  \thanks{This work was supported by ETH Zurich funds, the SNF AP Energy Grant \#160573, and the SFOE grant \#SI/501708 UNICORN.}

  \maketitle
  \thispagestyle{empty}
  \pagestyle{empty}

  \begin{abstract}
    In this paper, we study the stability and convergence of continuous-time Lagrangian saddle flows to solutions of a convex constrained optimization problem.
    Convergence of these flows is well-known when the underlying saddle function is either strictly convex in the primal or strictly concave in the dual variables. In this paper, we show convergence under non-strict convexity when a simple, unilateral augmentation term is added. For this purpose, we establish a novel, non-trivial characterization of the limit set of saddle-flow trajectories that allows us to preclude limit cycles.
    With our presentation we try to unify several existing problem formulations as a projected dynamical system that allows projection of both the primal and dual variables, thus complementing results available in the recent literature.
  \end{abstract}


\end{frontmatter}

\section{Introduction}\label{sec:intro}

Saddle flows have historically been a core element in many different domains, from circuit theory~\citep{smale1972mathematical, brayton1964theory1} to port-Hamiltonian systems~\citep{vanderschaft_port-hamiltonian_2014}, and their study dates back to the seminal works by~\cite{arrow_studies_1958, kose1956solutions, venets1985continuous}.

Currently, saddle flows are also used in the context of \emph{feedback} (or \emph{autonomous}) optimization, with
applications in power systems, network optimizations, etc. (see \citealt{colombino2019online, dallanese2018, hauswirth2019timescale, paganini2010}).
The goal of feedback optimization is to steer a physical system to an optimal operation point while satisfying operational constraints. Some of these need to be satisfied by the system trajectory at all times, while others need to be satisfied at steady-state.
Projections on the primal variables are needed to enforce trajectory constraints.
On the other hand, steady-state constraints often come as inequality constraints that can be dualized.
This leads to dual variables which have to be positive and therefore projections on the dual variables of the saddle flow are needed.

In this paper, we analyze the convergence and stability of projected saddle flows as they are encountered in feedback optimization. By considering projections on both the primal and the dual variables, we capture effects of physical saturation and discontinuity in the optimization strategy in a single mathematical model.

The contributions of this paper are as follows:
On the technical side, we expand on existing work and characterize the zero-dissipation set of a saddle flow. This allows us to characterize the emergence of limit cycles under non-strict convexity and show that these can be avoided with an augmentation term that does not affect the equilibria of the saddle flow. From a conceptual perspective, this shows how saddle flows can solve pure \emph{feasibility problems}.
On the educational side, we strive for a concise and accessible presentation and derive results that subsume other recent results for convex problems. In addition, we also hint at extensions of projected saddle flows on non-convex domains.

Similar to related work, our stability and convergence results are based on elementary methods such as basic convex analysis and a LaSalle invariance argument applied to a squared distance function. Beyond that, our new characterization of the zero-dissipation set requires additional facts about convexity and careful logical reasoning.
The only advanced technical aspect is the use of projected dynamical systems to model the discontinuous dynamics that arise in order to enforce constraints on primal and dual variables. In contrast to other frameworks that have been exploited for this purpose, such as monotone mappings~\citep{goebel_stability_2017}, complementarity systems~\citep{stegink2018convergence}, or subgradient formulations, existence and uniqueness of trajectories for projected dynamical systems is guaranteed for non-convex domains and even on abstract manifolds \citep{hauswirth_projected_2018-1}. We argue that this makes them a natural choice for future research on projected saddle flows that ventures beyond convexity.

Our paper differs from related work as follows:
\cite{cherukuri_2016_asymptotic} and \cite{cherukur_2018_role} do not have primal projections. \cite{niederlaender_distributed_2016} and~\citep{stegink2018convergence} do not have dual projections. \cite{dhingra2018proximal} works with non-convex and non-differentiable cost functions, but does not have any projections. These papers therefore do not analyze the same problem.
In~\cite{tangRunningPrimalDualGradient2018a}, the Lagrangian is regularized with an \emph{a priori} guess of the optimal dual variable. This makes the saddle flow strictly concave in the dual variable. This improves the convergence, but the resulting equilibrium is not the solution of the original optimization problem.
In contrast to~\cite{goebel_stability_2017} we assume that all relevant problem components are differentiable. This allows us to characterize the limit behavior in the case where the saddle-function is not strictly concave/convex in the dual nor the primal variable. This way, we are able to extend results in \cite{goebel_stability_2017} and \cite{Cherukuri2015}. The assumption that all relevant problem components are differentiable is motivated by the fact that feedback optimization controllers based on subgradient formulations are not easily implementable \citep{hauswirth2019timescale}.
The saddle flow limit behavior that we identify corresponds to what has also been predicted in \cite{holding_2014_oscillations}, where different methods have been used in the proof. The proposed augmentation of the Lagrangian to ensure convergence is also different.

The rest of the paper is organized as follows. In Section~\ref{sec:prelim} we fix the notation and review the basic properties of projected dynamical systems.
In Section~\ref{sec:prob_form} we provide the formal problem description. Sections~\ref{sec:stab} and~\ref{sec:zero-dis_set} discuss stability and convergence under strict convexity, and characterize the zero-dissipation set, respectively. In Section~\ref{sec:nonstrict} we present the limit behavior under non-strict convexity and in Section~\ref{sec:example} we give an example. Finally, we summarize our results and discuss open problems in Section~\ref{sec:conc}.

\section{Preliminaries}\label{sec:prelim}

\subsection{Notation}

We consider $\bbR^n$ with the usual inner product $\left\langle \cdot, \cdot \right\rangle$ and 2-norm $\| \cdot \|$. Hence, a linear map is given by a matrix $\bbR^{n \times m}$. We denote the \emph{closure} of a set $\calX \subseteq \bbR^n$ by $\cl \calX$.
Given a differentiable map $G: \bbR^n \rightarrow \bbR^m$ the \emph{Jacobian of $G$ at $x$} is the $n\times m$-matrix of partial derivatives denoted by $\nabla G(x)$. Namely, if $G: \bbR^n \rightarrow \bbR$ the \emph{gradient of $G$ at $x$} is $\nabla G(x)$. The partial Jacobian with respect to an argument or variable $y$ is denoted by $\nabla_y G$.
For a vector $v \in \bbR^n$ and an index set $I \subseteq \{1, \ldots, n\}$ we denote by $v_I$ the vector obtained by stacking the $I$-th components of $v$. Similarly for a matrix~$A$, $A_I$ denotes the matrix made up of the $I$-th columns of~$A$.
The \emph{Lie derivative of $G$} along vector field $f:\bbR^n \rightarrow \bbR^n$ at $x$, denoted by $\calL_f G(x)$ is the directional derivative of $G$ at $x$ in the direction of $f$. That is, $\calL_f G(x) := \nabla G(x)^T f(x)$. By $\bbR^m_{\geq0}$ we denote vectors of size $m$ with positive entries including 0.

\subsection{Projected Dynamical Systems}\label{ssec:proj_dyn_sys}

We quickly recall the definition and properties of projected dynamical systems on convex domains.

\begin{definition}\label{def:tgt_cone}
  Given a closed convex set $\calX \subseteq \bbR^n$, the \emph{tangent} and \emph{normal cone} at $x \in \calX$ are respectively defined as
  \begin{align}
    T_x \calX & := \cl \{ v \in \bbR^n \, | \,\exists \lambda > 0: \, x + \lambda v \in \calX \}                                                    \\
    N_x \calX & := \{ \eta \in \bbR^n \, | \, \forall v \in T_x \calX: \, \left\langle v, \eta \right\rangle \leq 0 \} \, .\label{eq:norm_cone_def}
  \end{align}
  The cones $T_x \calX$ and $N_x \calX$ are both closed, convex, and they are polar to each other.
\end{definition}

\begin{remark}
  The tangent cone takes an explicit form for sets $\calX := \{ x \, | \, \zeta(x) \leq 0 \}$ where $\zeta: \bbR^n \rightarrow \bbR^p$ and where $\calX$ satisfies constraint qualifications~\citep[Thm 6.14]{rockafellar_variational_1998}. Namely, in this case we have
  \begin{align*}
    T_x \calX & := \{ v \in \bbR^n \, | \, \nabla \zeta_{I(x)}(x)^T v \leq 0 \} \, ,
  \end{align*}
  where $I(x) := \{ i \, | \, \zeta_i(x) = 0 \}$ denotes the set of active constraints.
\end{remark}

We define a differential projection operator for a closed convex set $\calX \subseteq\bbR^n$, $x \in \calX$, and $v \in \bbR^n$ as
\begin{align*}
  \left[ v \right]_\calX^x :=  \arg \underset{w \in T_x \calX}{\min} \, \| v - w \|^2 \, ,
\end{align*}
that is, $[ v ]_\calX^x$ projects a vector $v$ onto the tangent cone of $\calX$ at the point $x$.
Since $T_x \calX$ is a closed convex set for any $x \in \calX$, the minimum norm projection of $v$ on $T_x \calX$ exists and is unique, and $[v]_{\calX}^x$ is well-defined. Furthermore, it holds that $\epsilon [v]_{\calX}^x = [ \epsilon v]_{\calX}^x$ for all $\epsilon > 0$ since $T_x \calX$ is a cone. Further, the following property is a consequence of Moreau's Theorem (e.g.,~\citep[Ex 12.22]{rockafellar_variational_1998}) and has been exploited for projected dynamical systems in~\cite{Cornet1983,Aubin1984,heemels_projected_2000, hauswirth_projected_2018-1} and others.

\begin{lemma}\label{lem:proj_decomp}
  For a closed convex set $\calX \subseteq \bbR^n$, $x \in \calX$, and $v \in \bbR^n$, it holds that $v - [v]_\calX^x \in N_x \calX$.
\end{lemma}

For a continuous vector field $F : \calX \rightarrow \bbR^n$, we define the projected dynamical system by applying the projection operator to the vector field $F$ at every point.
This leads to the initial value problem
\begin{align}\label{eq:pds_def}
  \dot x = [F(x)]_\calX^x \, , \quad x(0) = x_0 \, ,
\end{align}
where $x_0 \in \calX$ denotes an initial condition.

In general, $[ F(x) ]_\calX^x$ is not continuous and standard existence results for ordinary differential equations do not apply.
Instead, a (complete) \emph{(Carath\'eodory) solution} to~\eqref{eq:pds_def} is defined as an absolutely continuous function $x: [0, \infty) \rightarrow \calX$ and $x(0) = x_0$, and for which $\dot x(t) = [F(x(t))]_\calX^x$ holds almost everywhere, i.e., for almost all $t \in [0, \infty)$.
Note that a solution to~\eqref{eq:pds_def} has to be \emph{viable}, i.e., remain in $\calX$ for all $t \in [0, \infty)$ by definition.

Historically, the earliest existence results for projected dynamical systems appear to date back to~\cite{henry_existence_1973} (for convex domains) and~\cite{Cornet1983} (non-convex, but tangentially regular sets). Projected dynamical systems also appear as special cases in the study of differential inclusions~\citep{Aubin1984} and viability theory~\citep{aubinViabilitytheory1991}. For more recent works that consider projected dynamical systems on manifolds and Banach spaces see~\cite{hauswirth_projected_2018-1, daniele_panoramic_2010} and references therein.
In contrast, the independent line of research in~\cite{Nagurney1996}, popular in the context of variational inequalities, uses a different approach which cannot recover the generality of the earlier works.\footnote{In particular, the often-cited existence result~\citep[Thm 2.5]{Nagurney1996} requires~\citep[Thm 2.7]{Nagurney1996} which, strictly speaking, only applies to convex polyhedra.}

For our purpose, we state the following existence and invariance result that can be recovered as a special case and combination of more general results in~\cite{Aubin1984} and \cite{hauswirth_projected_2018-1} that do not require $\calX$ to be convex and $F$ Lipschitz:

\begin{theorem}\label{thm:lasalle}
  Consider~\eqref{eq:pds_def}. Let $\calX$ be closed convex, $F(x)$ locally Lipschitz, $\Psi: \bbR^{n} \rightarrow \bbR$ continuously differentiable, and $\mathcal{S}_\ell := \{ x \in \calX \, | \, \Psi(x) \leq \ell \}$ compact for all $\ell$. If $\calL_{\eqref{eq:pds_def}} \Psi(x) \leq 0$ for all $x \in \calX$, then there exists a unique complete solution of~\eqref{eq:pds_def} for every $x_0 \in \calS_\ell$ which converges to the largest invariant subset of $\cl \lbrace x \in \calS_\ell \, | \, \calL_{\eqref{eq:pds_def}} \Psi(x) = 0\rbrace$.
\end{theorem}

Note that $\Psi$ needs to have compact sublevel sets with respect to $\calX$, which is satisfied if $\calX$ is itself compact.

\begin{remark}
  Theorem~\ref{thm:lasalle} combines several fundamental results: First, existence of local solutions is guaranteed by standard viability results (e.g.~\citealt[Thm 3.3.4]{aubinViabilitytheory1991}). For an invariance principle such as~\citep[Thm 3]{bacciotti_stability_1999} to be applicable, \emph{continuous dependence on initial conditions} is required. This is the case because~\eqref{eq:pds_def} is equivalent to \emph{well-posed} differential inclusion (see \citealt[Thm 3, \S8]{filippovDifferentialEquationsDiscontinuous1988} and~\citealt[Thm 6.3]{hauswirth_projected_2018-1}). Finally, invariance of the compact sublevel sets of $\Psi$ guarantees the existence of complete solutions.
\end{remark}

\section{Problem Formulation}\label{sec:prob_form}

For the remainder of the paper, we will refer to the convex optimization problem given by
\begin{subequations}\label{eq:basic_opt}
  \begin{align}
    \underset{x}{\text{minimize}} \quad & f(x) \label{eq:opt_obj}              \\
    \text{subject to} \quad             & g(x) \leq 0 \label{eq:opt_soft}      \\
                                        & x \in \calX \, . \label{eq:opt_hard}
  \end{align}
\end{subequations}
Here, we make the following assumption.
\begin{assumption}\label{ass:convex_diff}
  The functions $f: \bbR^{n} \rightarrow \bbR$ and $g: \bbR^n \rightarrow \bbR^m$ are convex and differentiable, and $\calX \subseteq \bbR^n$ is closed and  convex.
\end{assumption}

In our context, $\calX$ defines a set of constraints that have to be satisfied at all times during the evolution of the dynamical system. This is important when a physical system is actuated in closed-loop with the saddle-flow algorithm.
In comparison, $g(x) \leq 0$ only has to be satisfied at steady-state.

\begin{assumption}[Well-posedness]\label{ass:bare_minimal}
  A (bounded) optimizer to~\eqref{eq:basic_opt} exists and Slater's condition holds, i.e., the feasible set $\{x \in \calX \, | \, g(x) \leq 0 \}$ has a non-empty relative interior.
\end{assumption}

Assumption~\ref{ass:bare_minimal} is standard to guarantee the applicability of the KKT conditions to certify global optimality~\citep{beck_introduction_2014}:

\begin{proposition}\label{prop:kkt_cond}
  Under Assumptions~\ref{ass:convex_diff} and~\ref{ass:bare_minimal}, a feasible point $x^\star$ is a global optimizer of~\eqref{eq:basic_opt} if and only if there exists $\mu^\star \in \bbR^m_{\geq0}$ such that
  \begin{align}\label{eq:kkt_stat}
    \nabla f(x^\star) + \nabla g(x^\star) \mu^\star \in - N_{x^\star} \calX
  \end{align}
  and for all $i=1,\ldots,m$
  \begin{align}\label{eq:kkt_comp}
    \mu^\star_i > 0 \quad \Rightarrow \quad g_i(x^\star) = 0 \, .
  \end{align}
\end{proposition}
Throughout the paper we refer to a primal-dual pair $(x^\star, \mu^\star)$ that satisfies the conditions of Proposition~\ref{prop:kkt_cond} as a \emph{solution} of~\eqref{eq:basic_opt}, whereas an~\emph{optimizer} of~\eqref{eq:basic_opt} only refers to the primal point $x^\star$. In particular, even if~\eqref{eq:basic_opt} has a unique optimizer (e.g., due to strict convexity of $f$) it may have multiple solutions because $\mu^\star$ might be not unique.

\begin{remark}
  We consider only inequality constraints $g(x) \leq 0$, since they introduce a discontinuous projection on the dual variables, which poses one of the main technical challenges. The case of equality constraints can be treated analogously with several simplifications (e.g., convex equality constraints are necessarily linear-affine and~\eqref{eq:kkt_comp} is vacuous) and has been well-documented in the literature.
\end{remark}

We now formally define the class of double-projection saddle flows that \emph{solve}~\eqref{eq:basic_opt}, i.e., whose trajectories converge to the global solutions of~\eqref{eq:basic_opt}, while satisfying $x \in \mathcal{X}$ at all times.

For this we define the \emph{state space} $\calZ := \calX \times \bbR^{m}_{\geq0}$ and the \emph{(primal-dual) feasible set}
\begin{equation}\label{eq:primal_dual_feasible_set}
  \calF := \{ (x, \mu) \in \calZ \, | \, g(x) \leq 0 \}.
\end{equation}
If $\calX$ is known to be a closed convex set, then it is immediate that $\calZ$ is a closed convex set with tangent cone $T_z \calZ = T_x \calX \times T_\mu \bbR^m_{\geq0}$ where $z := (x, \mu)$. Hence, we can define a well-behaved projected dynamical system restricted to $\calZ$.
For this, we first dualize the inequality constraint $g(x)$ in \eqref{eq:opt_soft}, but not the constraint in \eqref{eq:opt_hard} leading to a partial Lagrangian. We then augment this partial Lagrangian to get the augmented partial Lagrangian of~\eqref{eq:basic_opt} which is defined as
\begin{align*}
  L: \qquad \calZ & \rightarrow \bbR                                                           \\
  (x, \mu)        & \mapsto f(x) + \mu^T g(x) + \tfrac{\rho}{2} \| \max\{ 0, g(x) \} \|^2 \, ,
\end{align*}
where $\rho \geq 0$ is an augmentation parameter.

This Lagrangian is differentiable with derivative
\begin{align*}
  \begin{bmatrix} \nabla_x L \\ \nabla_\mu L \end{bmatrix}
  = \begin{bmatrix}
    \nabla f(x) + {\nabla g (x)}
    \left(
    \mu + \rho \max\{ 0, g (x) \}
    \right) \\
    g(x)
  \end{bmatrix} \, .
\end{align*}
We thus study the system described by projected gradient descent on the primal variable $x$ and a projected gradient ascent on the dual variable $\mu$, i.e.,
\begin{subequations} \label{eq:ahu_dyn}
  \begin{align}
    \dot x   & = \left[ - \nabla_x L(x, \mu) \right]^{x}_{\calX} \label{eq:ahu_dyn_a}                    \\
    \dot \mu & =    \left[ \nabla_\mu L(x, \mu) \right]^{\mu}_{\bbR^m_{\geq0}} \, . \label{eq:ahu_dyn_b}
  \end{align}
\end{subequations}
Existence of solutions for~\eqref{eq:ahu_dyn} will be a by-product of the forthcoming stability analysis which uses Theorem~\ref{thm:lasalle}. For now, note that $\nabla L$ is locally Lipschitz\footnote{Notice that $\max\{0, g(x)\}$ is convex and hence locally Lipschitz.} and $\calZ$ is closed convex, thus satisfying the requirements on $\calX$ and $F$ in Theorem~\ref{thm:lasalle}.

\begin{remark}
  Convex-concave saddle flows, such as~\eqref{eq:ahu_dyn}, have been studied using a variety of tools. In particular, the forthcoming results in Section~\ref{sec:stab} that exploit monotonicity to establish convergence to an invariant set, can be presented using complementarity systems~\citep{stegink2018convergence} or maximal monotone mappings~\citep{goebel_stability_2017}. The latter, in particular, allows for a very clean presentation, even if the saddle function $L$ is non-differentiable. However, when considering a non-convex problem (which is outside the scope of this paper) the monotonicity fails to hold. This not only jeopardizes stability, but also calls into question the applicability of existence results such as~\citep[Thm 1]{brogliatoequivalencecomplementaritysystems2006} or \citep[Thm 2.2]{goebel_stability_2017}.
\end{remark}

\begin{remark}
  The system~\eqref{eq:ahu_dyn} admits various variations that do not affect the main results presented in this paper. Namely,~\eqref{eq:ahu_dyn_a} and~\eqref{eq:ahu_dyn_b} can be subject to different time constants. Furthermore, instead of using a single augmentation parameter $\rho$, one can scale the augmentation term by a matrix $\Xi \succ 0$ as $\| \max\{0, g(x) \} \|^2_\Xi$ as in~\cite{stegink2018convergence}.
\end{remark}

For simplicity of notation we henceforth use
\begin{equation}\label{eq:def_vf}
  F(z) := \left[ \begin{smallmatrix} -\bbI_{n} & 0 \\ 0 & \bbI_{m} \end{smallmatrix} \right] \nabla L(z) \, .
\end{equation}
This definition alludes to the fact that~\eqref{eq:ahu_dyn} describes a projected \emph{pseudo-gradient flow} on $\calZ$, i.e., a gradient flow in an indefinite metric~\citep{bloch1992geometry, vanderschaft_port-hamiltonian_2014, smale1972mathematical}. Hence, we can write~\eqref{eq:ahu_dyn} compactly as
\begin{align}\label{eq:ahu_compact}
  \dot z = \left[ F(z) \right]_\calZ^z \, .
\end{align}
The following statement is well-known and can be easily derived from Proposition~\ref{prop:kkt_cond} and Lemma~\ref{lem:proj_decomp}:
\begin{proposition}\label{prop:equi_opti}
  Under Assumptions~\ref{ass:convex_diff} and~\ref{ass:bare_minimal} the set of equilibria of~\eqref{eq:ahu_compact} is equivalent to the set of global solutions of~\eqref{eq:basic_opt}.
\end{proposition}

\section{Stability of Projected saddle Flows}\label{sec:stab}

Next, we review the stability (but not necessarily convergence) of the dynamics~\eqref{eq:ahu_compact}. The results of this section use well-established tools and have been shown, in one way or another, in~\cite{cherukuri_2016_asymptotic,cherukur_2018_role,niederlaender_distributed_2016,stegink2018convergence,Cherukuri2015} and exploit basic facts from convex analysis~\citep{beck_introduction_2014, rockafellar_variational_1998}. We aim for a particularly concise presentation that prepares for the main result in the next section. For completeness, the proofs for this section are included in the appendix.

As a first step, it suffices to characterize the \emph{monotonicity} of the ``unprojected'' vector field $F(z)$ by exploiting the convexity of $g$ and $f$. The first two points in the following proposition are well-known, the last point is new (to the best of the authors' knowledge) but easy to prove and crucial for our forthcoming results. For completeness the proof can be found in the appendix.

\begin{proposition}\label{prop:semi-monot}
  Under Assumptions~\ref{ass:convex_diff} and~\ref{ass:bare_minimal} the following holds for $F(z)$ defined in~\eqref{eq:def_vf}:
  \begin{enumerate}[label=(\roman*)]
    \item\label{it:semi_monot_1} For all $z, \hat{z} \in \calZ$ the vector field $F(z)$ satisfies
          \begin{align*}
            \left\langle z - \hat{z}, F(z) - F(\hat{z}) \right\rangle \leq 0 \, .
          \end{align*}
    \item\label{it:semi_monot_2} If $f$ is strictly convex, for $z, \hat{z} \in \calZ$ with $x \neq \hat{x}$ one has
          \begin{align*}
            \left\langle z - \hat{z}, F(z) - F(\hat{z}) \right\rangle < 0 \, .
          \end{align*}
    \item\label{it:semi_monot_3} If $\rho > 0$, then for all $z \in \calZ$ but $z \notin \calF$ and all $\hat{z} \in \calF$ we have    \
          \begin{align*}
            \left\langle z - \hat{z}, F(z) - F(\hat{z}) \right\rangle < 0 \, .
          \end{align*}
  \end{enumerate}
\end{proposition}

One useful property of projected dynamical systems on convex sets is the fact that monotonicity of vector fields is preserved when passing from an unprojected to a projected vector field as the following lemma shows. For completeness the proof can be found in the appendix.

\begin{lemma}\label{lem:preserv_monot}
  Let $\calZ \subseteq \bbR^n$ be a closed convex set and $F: \calZ \rightarrow \bbR^n$ a vector field. If for $z, \hat{z} \in \calZ$ it holds that
  \begin{align*}
    \left\langle z - \hat{z}, F(z) - F(\hat{z}) \right\rangle \leq \alpha
  \end{align*}
  for some $\alpha \in \bbR$, then it also holds that
  \begin{align*}
    \left\langle z - \hat{z}, [F(z)]_\calZ^z - [F(\hat{z})]_\calZ^{\hat{z}} \right\rangle \leq \alpha \, .
  \end{align*}
\end{lemma}
Lemma~\ref{lem:preserv_monot} allows to amend Proposition~\ref{prop:semi-monot} as follows.
\begin{corollary}\label{prop:semi-monot_proj}
  The statements~\ref{it:semi_monot_1},~\ref{it:semi_monot_2}, and~\ref{it:semi_monot_3} in Proposition~\ref{prop:semi-monot} are valid if $F(z)$ and $F(\hat{z})$ are replaced with $[F(z) ]_\calZ^z$ and $[F(\hat{z}) ]_\calZ^{\hat{z}}$ respectively, everywhere.
\end{corollary}

Hence, the invariance principle for projected dynamical systems (Theorem~\ref{thm:lasalle}) can be used to conclude the following. For completeness the proof can be found in the appendix.

\begin{theorem}\label{thm:lasalle_appl}
  Under Assumptions~\ref{ass:convex_diff} and~\ref{ass:bare_minimal} there exists a unique complete solution $z:[0, \infty) \rightarrow \calZ$ of~\eqref{eq:ahu_compact} for every initial condition and $z$ converges asymptotically to the largest invariant subset $\Omega$ of $\calM_{z^\star}$ where
  \begin{align}\label{eq:zero-dissip-set}
    \calM_{z^\star} :=  \cl \{ z \in \calZ \, | \, \left\langle z - z^\star, [ F(z) ]_{\calZ}^z \right\rangle = 0 \} \, ,
  \end{align}
  where $z^\star$ is a solution of~\eqref{eq:basic_opt}.
  Moreover, all solutions $z^\star$ are stable.
\end{theorem}

If $f$ is strictly convex and with the help of item~(ii) of Proposition~\ref{prop:semi-monot} and Corollary~\ref{prop:semi-monot_proj} one can show convergence to the solutions of~\eqref{eq:basic_opt}.
This results was established for a more general setup in~\cite{goebel_stability_2017} for saddle-functions that are strictly convex/concave in either the primal or dual variables. In our case, the saddle-function, which is a partial Lagrangian of~\eqref{eq:basic_opt}, is never strictly concave in the dual variables. Hence, strict convexity has to stem from the primal variables.

\begin{corollary}\label{corr:asymp_cvg}\cite[Thm 3.3] {goebel_stability_2017}
  Under Assumptions~\ref{ass:convex_diff} and~\ref{ass:bare_minimal} and if $f$ is strictly convex, every trajectory of~\eqref{eq:ahu_compact} converges to a primal-dual solution of~\eqref{eq:basic_opt}.
\end{corollary}

\begin{remark}
  It is worth pointing out a potential pitfall: After concluding that trajectories converge to an invariant subset $\Omega$ of $\calM_{z^\star}$, it is tempting to study the dynamics restricted to $\Omega$, e.g., with another LaSalle function, and conclude convergence of any trajectory to this ``nested'' $\omega$-limit set. This approach is however in general not valid unless additional assumptions are satisfied~\citep{arsie_locating_2010}. Second, although the zero-dissipation set $\calM_{z^\star}$ (for strictly convex $f$) consists only of equilibrium points, Theorem~\ref{thm:lasalle_appl} does not imply convergence to a single point. To conclude pointwise convergence, one needs to exploit the definitions of stability and limit sets as done in~\cite{cherukuri_2016_asymptotic, goebel_stability_2017, stegink2018convergence}.
\end{remark}

\section{Characterization of the Zero-Dissipation Set}
\label{sec:zero-dis_set}
In order to better describe the limit behavior of~\eqref{eq:ahu_compact}, we now introduce a novel characterization of the zero-dissipation set. The following assumption will be required below for item~\ref{it:non-dissip_7} of Proposition~\ref{prop:zero-dissip_charact}. Whether it can be relaxed remains an open question.

\begin{assumption}[Strict complementary slackness]\label{ass:strict_comp}
  \leavevmode \\
  Under Assumption~\ref{ass:bare_minimal} there is a solution $(x^\star, \mu^\star)$ of \eqref{eq:basic_opt} such that
  \begin{align}\label{eq:kkt_strict_comp}
    \forall i=1,\ldots,m: \quad g_i(x^\star) = 0 \quad \Rightarrow \quad \mu_i^\star > 0 \, .
  \end{align}
\end{assumption}
Assumption~\ref{ass:strict_comp} is weak, because it only needs to hold for a single solution of~\eqref{eq:basic_opt}. Further, for large classes of (parametrized) optimization problems this assumption holds generically, i.e., for almost all problem instances~\citep{spingarn1979generic}.

We now recall a basic, but little known, result about convex functions which is crucial for our analysis.
\begin{lemma}~\citep[Ex. 7.28iii]{beck_introduction_2014}\label{lem:cvx_grad_eq}
  If $f: \bbR^n \rightarrow \bbR$ is convex and differentiable with Lipschitz gradient, then it holds $\forall \, x,y \in \bbR^n$ that
  \begin{align*}
    f(y) - f(x) = \nabla f(x)^T(x-y)  \quad \Rightarrow \quad \nabla f(x) = \nabla f(y)
  \end{align*}
  as well as
  \begin{align*}
    \left\langle x - y, \nabla f(x) - \nabla f(y) \right\rangle = 0 \quad \Rightarrow \quad \nabla f(x) = \nabla f(y) \, .
  \end{align*}
\end{lemma}
Hence, our key technical result reads as follows.
\begin{proposition}[Characterization of zero-dissipation set]\label{prop:zero-dissip_charact}
  \leavevmode \\
  Let Assumptions~\ref{ass:convex_diff} and~\ref{ass:bare_minimal} hold, let $z^\star=(x^\star, \mu^\star)$ be any solution of~\eqref{eq:basic_opt} and let $\calF$ be defined as in \eqref{eq:primal_dual_feasible_set}.
  Then the following statements about the zero-dissipation set $\calM_{z^\star}$ defined in~\eqref{eq:zero-dissip-set} hold true:
  \begin{enumerate}[label=(\roman*)]
    \item\label{it:non-dissip_6} If $\rho > 0$ or $f$ is strictly convex, then $\calM_{z^\star} \subseteq \calF$.
    \item\label{it:non-dissip_1} For all $(x, \mu) \in \calM_{z^\star}$ it holds that $\nabla f(x) = \nabla f(x^\star)$.
    \item\label{it:non-dissip_2} For all $(x, \mu) \in \calM_{z^\star}$ and all $i = 1, \ldots, m$ we have
          \begin{align*}
            {\mu^\star_i} = 0 \quad \text{or} \quad g_i(x) = g_i({x}^\star) + \nabla g_i({x}^\star) (x - {x}^\star) \, .
          \end{align*}
    \item\label{it:non-dissip_3} For all $(x, \mu) \in \calM_{z^\star}$ it holds that $\nabla g(x) \mu = \nabla g(x^\star) \mu$.
    \item\label{it:non-dissip_5} For all $(x, \mu) \in \calM_{z^\star}$ and all $i = 1, \ldots, m$ we have
          \begin{align*}
            g_i(x^\star) < 0  \quad \Rightarrow \quad \mu_i = 0 \, .
          \end{align*}
    \item\label{it:non-dissip_7}  If, in addition, Assumption~\ref{ass:strict_comp} is satisfied and the zero-dissipation set $\calM_{z^\star}$ in~\eqref{eq:zero-dissip-set} is defined with respect to a solution $z^\star$ that satisfies~\eqref{eq:kkt_strict_comp}, then it holds that for all $(x, \mu) \in \calM_{z^\star}$ and all $i = 1, \ldots, m$ we have
          \begin{align*}
            g_i(x^\star) = 0  \quad \Rightarrow \quad \mu_i > 0 \quad \text{or} \quad g_i(x) \geq 0\, .
          \end{align*}
  \end{enumerate}
\end{proposition}
\begin{pf}
  We use the definitions of $M_1, M_2, M_3$ from the proof of Proposition~\ref{prop:semi-monot} and the definitions of $\eta, \hat{\eta}$ from the proof of Lemma~\ref{lem:preserv_monot}, which can be found in the appendix.
  Based on \eqref{eq:zero-dissip-set}, we consider the case in which $\hat z = z^\star$ and $\hat\eta= \eta^\star$.

  For zero dissipation we need $M_1 = M_2 = M_3 = 0$ and
  \begin{align}\label{eq:proj_zero_dissip}
    \left\langle z - \hat{z}, \eta \right\rangle =\left\langle z - \hat{z}, \hat{\eta} \right\rangle = 0 \, .
  \end{align}
  First, if $\rho > 0$, then $M_2 = 0$ is only the case when $g(x) \leq 0$ and hence $z \in \calF$. Otherwise, if $f$ is strictly convex the primal optimizer of~\eqref{eq:basic_opt} is unique and $M_1 > 0$ for all $x \neq \hat{x}$ and in particular for $x \notin \calF$. This establishes~\ref{it:non-dissip_6}.

  Second, for $M_1 = \left\langle x - \hat{x}, \nabla f(x) - \nabla f(\hat{x}) \right\rangle = 0$ we can use Lemma~\ref{lem:cvx_grad_eq} to conclude~\ref{it:non-dissip_1}.

  Third, for $M_3 = 0$ to hold we require that
  \begin{align}
    \mu^T \left( - g(\hat{x}) + g(x)  + \nabla g(x)^T ( \hat{x} - x) \right)            & = 0 \label{eq:m3_zero_dissip1}      \\
    \hat{\mu}^T \left( - g(x) + g(\hat{x}) + \nabla g(\hat{x})^T (x - \hat{x} ) \right) & = 0 \, . \label{eq:m3_zero_dissip2}
  \end{align}
  From~\eqref{eq:m3_zero_dissip2} we directly infer~\ref{it:non-dissip_2}.

  From~\eqref{eq:m3_zero_dissip1}, we establish, using Lemma~\ref{lem:cvx_grad_eq}, that
  \begin{align*}
    \mu_i = 0 \quad & \text{or} \quad \nabla g_i(x) = \nabla g_i(\hat{x}) \, ,
  \end{align*}
  which implies~\ref{it:non-dissip_3}.

  Next, for~\ref{it:non-dissip_5}, recall that $\calZ := \calX \times \bbR^m_{\geq 0}$. Therefore, the projection on $\calZ$ can be decomposed into a projection on $\calX$ and $m$ projections on $\bbR_{\geq 0}$. Therefore, $\left\langle z - \hat{z}, \eta \right\rangle = 0$ implies that $\left\langle \mu_i - \hat{\mu}_i, \eta_{\mu i} \right\rangle = 0$ for all $i = 1, \ldots, m$
  where $\eta_{\mu i}$ denotes the $i$-th dual components of $\eta, \hat{\eta}$, respectively. The analogous statement holds for $\hat{\eta}_{\mu i}$.

  Combining these statements, we have, for all $i = 1, \ldots, m$,
  \begin{align}\label{eq:raw_complementarity}
    \eta_{\mu i} \neq 0 \quad \text{or} \quad \eta_{\hat{\mu} i} \neq 0 \quad \Rightarrow \quad \mu_i = \hat{\mu}_i \, .
  \end{align}

  Note that for $\bbR_{\geq 0}$ we have $N_x \bbR_{\geq 0} = \bbR$ for all $x > 0$ and $N_0 \bbR_{\geq 0} = \bbR_{\geq 0}$. Thus, using Lemma~\ref{lem:proj_decomp},  we have that
  \begin{align*}
    \hat{\eta}_{\mu i} = \begin{cases}
      0            & \quad \text{if } g_i(\hat{x}) \geq 0 \text{ or } \hat{\mu}_i > 0 \\
      g_i(\hat{x}) & \quad \text{otherwise} \, .
    \end{cases}
  \end{align*}
  The statement~\ref{it:non-dissip_5} follows immediately because, since $\hat{x}$ is a solution of~\eqref{eq:basic_opt}, $g_i(\hat{x}) < 0$ implies $\hat{\mu}_i = 0$ by complementary slackness~\eqref{eq:kkt_comp} and therefore $\eta_{\mu i} \neq 0$, and with~\eqref{eq:raw_complementarity} it follows that $\mu_i = \hat{\mu}_i$.

  Finally, for~\ref{it:non-dissip_7} we work with the contraposition of~\eqref{eq:raw_complementarity} and the strict complementarity slackness assumption on $(\hat{x}, \hat{\mu})$. Namely, if $g_i(\hat{x}) = 0$ we have that $\hat{\eta}_{\mu i} = 0$ and $\hat{\mu}_i > 0$ by~\eqref{eq:kkt_strict_comp}. If $\mu_i = 0$, then~\eqref{eq:raw_complementarity} implies that $\eta_{\mu i} \neq 0$ which in turn implies that $g_i(x) \geq 0$.
  \qed
\end{pf}

Proposition~\ref{prop:zero-dissip_charact} is instrumental to study the occurrence of limit cycles for non-strictly convex cost functions.

\section{Limit Behavior under Non-strict Convexity}\label{sec:nonstrict}

When $f$ is not strictly convex it is well-known that convergence to equilibrium points is not in general guaranteed and limit cycles can occur. In the following, we do not consider projections of the primal variables and hence make the assumption that $\calX = \bbR^n$. Furthermore, we make use of Assumption~\ref{ass:strict_comp}.

With the help of Proposition~\ref{prop:zero-dissip_charact} we can show that on the zero dissipation set~$\calM_{z^\star}$ the dynamics~\eqref{eq:ahu_compact}  follow simple linear Hamiltonian dynamics.
Since we are only interested in invariant subsets of $\calM_{z^\star}$ we can then establish a contradiction which leads us to conclude that, in the presence of augmentation of the Lagrangian, the only invariant subsets of $\calM_{z^\star}$ are equilibrium points.

\begin{proposition}\label{prop:local_hamiltonian}
  Let Assumptions~\ref{ass:convex_diff},~\ref{ass:bare_minimal} and~\ref{ass:strict_comp} hold, let $\calX = \bbR^n$, and let the zero-dissipation set $\calM_{z^\star}$ in~\eqref{eq:zero-dissip-set} be defined with respect to a solution $z^\star$ that satisfies~\eqref{eq:kkt_strict_comp}. Then, for all $z \in \calM_{z^\star}$, the dynamical system~\eqref{eq:ahu_compact} reduces to
  \begin{align}\label{eq:local_lin_ode}
    \begin{bmatrix} \dot x \\ \dot \mu_{I^\star} \end{bmatrix} & = \begin{bmatrix} 0 & -A \\ A^T & 0 \end{bmatrix}
    \begin{bmatrix} x \\ \mu_{I^\star} \end{bmatrix} - \begin{bmatrix} c\\ d \end{bmatrix} \, ,
  \end{align}
  where $c = \nabla f(x^\star)$, $d = \nabla g_{I^\star}(x^\star)^T  x^\star$, and $A = \nabla g_{I^\star}(x^\star) $, with $I^\star := \{ i \, | \, g_i(x^\star) = 0 \}$ and $\mu_{i} = 0$ for all $i \notin I^\star$.
\end{proposition}
\begin{pf}
  First, note that \ref{it:non-dissip_5} in Proposition~\ref{prop:zero-dissip_charact} implies that $\mu_i = 0$ for all $i \notin I^\star$ on all of $\calM_{z^\star}$.
  Second, using~\ref{it:non-dissip_6},~\ref{it:non-dissip_1},~\ref{it:non-dissip_3}, and~\ref{it:non-dissip_5} in Proposition~\ref{prop:zero-dissip_charact} we know that for all $(x, \mu) \in \calM_{z^\star}$ we have
  \begin{align*}
    - \nabla_x L(x, \mu) & = -  (\nabla f(x) + \nabla g(x) \mu)                                     \\
                         & = - (\nabla f(x^\star) + \nabla g(x^\star) \mu)                          \\
                         & = - (\nabla f(x^\star) + \nabla g_{I^\star}(x^\star) \mu_{I^\star}) \, .
  \end{align*}
  Furthermore, we know that $i \in I^\star$ implies $\mu^\star_i > 0$ according to~\eqref{eq:kkt_strict_comp}. Hence, we apply~\ref{it:non-dissip_2} in Proposition~\ref{prop:zero-dissip_charact} to conclude that for all $i \in I^\star$ we have
  \begin{align*}
    g_i(x) = g_i(x^\star) + \nabla g_i(x^\star)^T(x - x^\star)  = \nabla g_i(x^\star)^T(x - x^\star)\, ,
  \end{align*}
  that is, $g_i$ is linear on $\calM_{z^\star}$ and equal to its linearization at $x^\star$.

  Finally, with~\ref{it:non-dissip_7} of Proposition~\ref{prop:zero-dissip_charact} we know that $\mu_i > 0$ or $g_i(x) \geq 0$ holds on all of $\calM_{z^\star}$ and therefore the projection of the dual variables in $I^\star$ is never active, i.e.,
  \begin{align*}
    \dot \mu_i = [g_i(x)]_{+}^{\mu_i} = g_i(x) \quad
    \text{for }i \in I^\star
  \end{align*}
  and hence, on $\calM_{z^\star}$ the system~\eqref{eq:ahu_compact} reduces to~\eqref{eq:local_lin_ode}.
  \qed
\end{pf}
Note that Proposition~\ref{prop:local_hamiltonian} states that trajectories, while in $\calM_{z^\star}$, satisfy~\eqref{eq:local_lin_ode}, but that does not imply that $\calM_{z^\star}$ is invariant. In this sense, Proposition~\ref{prop:local_hamiltonian} does not prove the existence of limit cycles. In fact, as we show next, in the presence of an augmentation term the only invariant subsets of $\calM_{z^\star}$ are equilibrium points.
\begin{theorem}
  Under Assumptions~\ref{ass:convex_diff},~\ref{ass:bare_minimal} and~\ref{ass:strict_comp}, if $\calX = \bbR^n$ and $\rho > 0$, then the largest invariant set $\Omega \subseteq \calM_{z^\star}$ in Theorem~\ref{thm:lasalle_appl} is equivalent to the set of equilibria of~\eqref{eq:ahu_compact}. Furthermore, all trajectories of~\eqref{eq:ahu_compact} converge asymptotically to the set of solutions of~\eqref{eq:basic_opt}.
\end{theorem}
\begin{pf}
  The linear system~\eqref{eq:local_lin_ode} has analytic solutions
  \begin{align*}
    x(t)     & := x^\star + V \begin{bmatrix}
      \diag(\beta) \sin(\sigma t + \phi) \\ \gamma \end{bmatrix}        \\
    \mu_I(t) & := \mu^\star - U \begin{bmatrix} \diag(\beta) \cos(\sigma t + \phi) \end{bmatrix} \, ,
  \end{align*}
  where $A$ has the singular value decomposition $A = U \Sigma V^T$ with $\Sigma = \begin{bmatrix} \diag(\sigma) & 0 \end{bmatrix}$, $\cos$ and $\sin$ apply componentwise, and $\beta, \phi \in \bbR^{|I^\star|}$, $\gamma \in \bbR^{n - |I^\star|}$ depend on initial conditions.

  For any initial condition and any $\beta \neq 0$ there exists $t$ such that for at least one component $i \in I^\star$ we have that $\dot \mu_i(t) >0$.
  Recall that, on $\calM_{z^\star}$,
  \begin{align*}
    \dot \mu_i(t) = A x - d = \nabla g_i(x^\star)(x - x^\star) = g_i(x),
  \end{align*}
  and therefore, unless $\beta = 0$, there exists at least one oscillating mode that drives the trajectory outside of the feasible set. However, by~\ref{it:non-dissip_6} in Proposition~\ref{prop:zero-dissip_charact} we have that $\calM_{z^\star} \subseteq \calF$ and therefore such an orbit is not invariant with respect to $\calM_{z^\star}$. Hence, the only invariant subsets of $\calM_{z^\star}$ are equilibrium points. Conversely, any equilibrium point is invariant and contained in $\calM_{z^\star}$. Since $\Omega$ is the \emph{largest} invariant subset of $\calM_{z^\star}$ it has to contain all equilibria.
  \qed
\end{pf}

\section{Illustrative Example}
\label{sec:example}
\begin{figure}
  \vspace{.3cm}
  \input{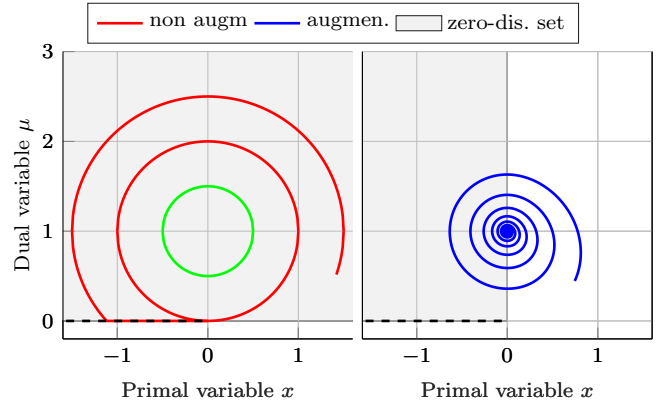}
  \caption{Left plot without augmentation ($\rho=0$): Trajectory (red) of \eqref{eq:example_dyn} converging to the boundary of the largest invariant subset $\Omega=\{(x,\,u)\in\bbR \times \bbR_{\geq0}\, |\, \|(x,\mu-1)\|_2<1\}$. Exemplary trajectory (green) of \eqref{eq:example_dyn} inside of $\Omega$. The zero-dissipative set $\calM_{z^\star}$ is everything besides the negative x-axis.
    Right plot with augmentation ($\rho>0$): Trajectory (blue) of \eqref{eq:example_dyn} converging to the largest invariant subset $\Omega=\{(x^\star,\mu^\star)\}$. The zero-dissipative set $\calM_{z^\star}$ is the second quadrant without the negative x-axis.}
  \label{fig:example}
\end{figure}

To illustrate the occurrence of limit cycles and that we do not experience such behavior with an augmented Lagrangian we analyze the stylized problem \[ \max_{x\le 0}\quad x \]
with the unique solution $(x^\star, \mu^\star)=(0,1)$.
We dualize the constraint and augment the Lagrangian.
Hence, the projected saddle flow~\eqref{eq:ahu_compact} for this problem is defined on $\calZ = \bbR \times \bbR_{\geq0}$ and given by
\begin{align}\label{eq:example_dyn}
  \dot{z}=\begin{bmatrix}
    \dot{x} \\
    \dot{\mu}
  \end{bmatrix}= \begin{bmatrix}
    1-\mu-\rho \max \{0,x\}\\{[x]}^{\mu}_{\bbR_{\geq0}}
  \end{bmatrix} \, ,
\end{align}
where $\rho \geq 0$ is the augmentation parameter.
A LaSalle function is given by
\begin{align*}
   & V(x,\mu) = \frac{1}{2}\|x-x^\star\|^2 + \frac{1}{2} \|\mu-\mu^\star\|^2 \geq 0 \, ,
\end{align*}
and the Lie derivative of $V(x,\mu)$ is
\begin{multline*}
  \calL_{\eqref{eq:example_dyn}} V(x,\mu)
  = (x - {[x]}^{\mu}_{\bbR_{\geq0}})(1 - \mu) - \rho x \max \{0,x\}\\
  =\begin{cases}
    x-\max\{0,x\}-\rho x \max \{0,x\} & \mu=0   \\
    -\rho x \max\{0,x\}               & \mu >0.
  \end{cases}
\end{multline*}
For $\rho=0$ (no augmentation) the zero-dissipation set is $\calM_{z^\star} = \{ z \in \calZ\} \setminus\{x<0,\mu=0\}$, and therefore there is only dissipation on the negative $x$-axis.
The red line in Fig.~\ref{fig:example} shows that in this case all trajectories converge to the circle around the solution $(0,1)$ with radius 1. On this set the trajectories form periodic orbits around the solution. With augmentation the zero-dissipation set is $\calM_{z^\star} = \{z \in \calZ \, | \,x\leq 0, \mu >0 \}$. Due to the dissipation on the negative $x$-axis and in the first quadrant, all trajectories converge to the largest invariant subset $\Omega$, which contains only the solution. Compare the blue line in Fig. \ref{fig:example}.

\section{Conclusion}\label{sec:conc}

In this paper, we proposed a unified formulation of saddle dynamics in the presence of projection of both the primal and the dual variables.

Such a unified approach, based on the formalism of projected dynamical systems, allows us to derive stability and convergence results that subsume more specific results that have been presented in the literature.

Also, we characterize the zero-dissipation set of a saddle flow and the emergence of limit cycles under non-strict convexity. We show that oscillations can be avoided with a simple augmentation term that does not affect the equilibrium of the saddle flow.

We expect to extend the analysis of these dynamics to non-convex domains, profiting from the well-posedness of projected dynamical systems on these domains.
This extension has the potential to support the application of these  methods for the analysis and design of feedback optimization schemes for a wide class of systems.

\bibliography{IEEEabrv,bibliography}

\appendix

\section{Proof of Proposition~\ref{prop:semi-monot}}
The proof mainly exploits the definition of convexity and follows ideas from~\cite{cherukuri_2016_asymptotic, stegink2018convergence} and others. Using the definition of $\nabla L$ we can decompose
\begin{align*}
  \left\langle z - \hat{z}, F(z) - F(\hat{z})) \right\rangle = - M_1 - M_2 - M_3 \, ,
\end{align*}
where we define
\begin{align*}
  M_1 & :=  \left\langle x - \hat{x}, \nabla f(x) - \nabla f(\hat{x}) \right\rangle                   \\
  M_2 & :=  \rho \left\langle x - \hat{x}, \nabla g(x) \max\{0, g(x)\}\right\rangle                   \\
      & \quad -  \rho \left\langle x - \hat{x}, \nabla g(\hat{x}) \max\{0, g(\hat{x})\} \right\rangle \\
  M_3 & := \left\langle x - \hat{x}, \nabla g(x) \mu - \nabla g(\hat{x}) \hat{\mu} \right\rangle      \\
      & \qquad - \left\langle \mu - \hat{\mu}, g(x) - g(\hat{x})\right\rangle \, .
\end{align*}
Hence, for the first term we have by the definition of convexity of $f$ that $M_1 \geq 0$, and if $f$ is strictly convex we have $M_1 > 0$ for all $x \neq \hat{x}$.

For $M_2$ we notice that if $\hat{z} \in \calF$ we have that $\max\{0, g(\hat{x})\} = 0$. Furthermore, we define $I(x):=\{ i \, | \, g_i(x)=0 \}$ and notice that if $z \in \calX$, but $z \notin \calF$, then for all $i \in I(x)$, where $I(x)\neq \emptyset$, we have $g_i(x)>0$. Therefore, using also convexity of $g$, it holds that
\begin{align*}
  M_2 & = \rho \left\langle x - \hat{x}, \nabla g(x) \max\{0, g(x)\} \right\rangle                                                      \\
      & = \rho \sum_{i\in I(x)}  (x_i - \hat{x}_i)^T \nabla g_i(x) g_i(x)                                                               \\
      & \geq  \rho \sum_{i\in I(x)}  (\underbrace{g_i(x)}_{>0} - \underbrace{g_i(\hat{x})}_{\leq 0}) \underbrace{g_i(x)}_{> 0} > 0 \, .
\end{align*}
Finally, for $M_3$ we have
\begin{align*}
  M_3 & = \mu^T \nabla g(x)^T (x - \hat{x}) -  \hat{\mu}^T \nabla g(\hat{x})^T (x - \hat{x})                                  \\ & \qquad \qquad
  - \mu^T (g(x) - g(\hat{x})) + \hat{\mu}^T(g(x) - g(\hat{x}))                                                                \\
      & = \mu^T \underbrace{\left( g(\hat{x}) - g(x)  - \nabla g(x)^T ( \hat{x} - x) \right)}_{\geq 0 \text{ (by convexity)}} \\ & \qquad \qquad
  + \hat{\mu}^T \underbrace{\left( g(x) - g(\hat{x}) - \nabla g(\hat{x})^T (x - \hat{x} ) \right)}_{\geq 0 \text{ (by convexity)}}\, ,
\end{align*}
and since $\mu, \hat{\mu} \geq 0$ it follows that $M_3 \geq 0$.
Hence,~\ref{it:semi_monot_1},~\ref{it:semi_monot_2}, and~\ref{it:semi_monot_3} follow immediately.

\section{Proof of Lemma~\ref{lem:preserv_monot}}

Lemma~\ref{lem:proj_decomp} states that there are normal vectors $\eta \in N_z \calZ$ and $\hat{\eta} \in N_{\hat{z}} \calZ$ such that $[F(z)]_{\calZ}^z = F(z) - \eta$ and $[F(\hat{z})]_{\calZ}^{\hat{z}} = F(\hat{z}) - \hat{\eta}$. Further, by definition of the normal cone to a convex set we have $\left\langle \eta, \hat{z} - z \right\rangle \leq 0$ and $\left\langle \hat{\eta}, z - \hat{z}\right\rangle \leq 0$ $\forall z,\hat{z} \in \calZ$. Hence,
\begin{align*}
   & \left\langle z - \hat{z}, [F(z)]_\calZ^z - [F(\hat{z})]_\calZ^{\hat{z}} \right\rangle                                                                             \\
   & \quad = \left\langle z - \hat{z}, F(z) - \eta - F(\hat{z}) + \hat{\eta} \right\rangle                                                                             \\
   & \quad = \underbrace{\left\langle z - \hat{z}, F(z) - F(\hat{z}) \right\rangle}_{\leq \alpha} + \underbrace{\left\langle \hat{z} - z, \eta \right\rangle}_{\leq 0}
  + \underbrace{\left\langle z - \hat{z}, \hat{\eta} \right\rangle}_{\leq 0}
\end{align*}
which immediately proves the lemma.

\section{Proof of Theorem~\ref{thm:lasalle_appl}}

By Assumption~\ref{ass:bare_minimal} there exists an optimizer $x^\star$ to~\eqref{eq:basic_opt} with dual solution $\mu^\star$ satisfying the conditions of Proposition~\ref{prop:kkt_cond} and $(x^\star, \mu^\star)$ is an equilibrium of~\eqref{eq:ahu_compact} by Proposition~\ref{prop:equi_opti}. Using Proposition~\ref{prop:semi-monot_proj} we know that for every (local) solution $z:[0, T) \rightarrow \calZ $ of~\eqref{eq:ahu_compact} we have for almost all $t \in [0, T)$ that
\begin{align*}
  \frac{d}{dt} \| z(t) - z^\star \|^2 & =
  2 \left\langle z - z^\star, [F(z)]_\calZ^z - [F(z^\star)]_\calZ^{z^\star} \right\rangle \\
                                      & =
  2 \left\langle z - z^\star, [F(z)]_\calZ^z  \right\rangle \leq 0 \, .
\end{align*}
Since $\| z - z^\star \|^2$ has compact level sets as a function of $z$ on $\calZ$ it follows from Theorem~\ref{thm:lasalle} that all trajectories are complete and converge to the largest invariant subset $\Omega$ of $\calM_{z^\star}$.
The non-positivity of the Lie derivative of $\| z(t) - z^\star \|^2$ guarantees stability via a standard Lyapunov argument.
\end{document}